\documentclass[10pt,twoside]{article}
\usepackage{amsmath,amscd}
\usepackage{amssymb}
\usepackage{epsfig}
\usepackage{Latex-document}

\def\volumeno{I}
\title{\bf Pattern Theory: The Mathematics of \vskip -2mm
Perception\vskip 6mm}
\author{David Mumford\vspace*{-0.5cm}\thanks{Division of Applied
Mathematics, Brown University, Providence RI 02912, USA. E-mail:
David\_Mumford@brown.edu}}
\date{\vspace{-8mm}}

\begin{document}

\maketitle

\thispagestyle{first} \setcounter{page}{401}

\begin{abstract}

\vskip 3mm

Is there a mathematical theory underlying intelligence? Control theory
addresses the output side, motor control, but the work of the last 30
years has made clear that perception is a matter of Bayesian statistical
inference, based on stochastic models of the signals delivered by our
senses and the structures in the world producing them. We will start by
sketching the simplest such model, the hidden Markov model for speech,
and then go on illustrate the complications, mathematical issues and
challenges that this has led to.

\vskip 4.5mm

%\noindent {\bf 2000 Mathematics Subject Classification:} 60K35, 82B43.

\end{abstract}

\vskip 12mm

\section{Introduction} \label{section 1}\setzero
\vskip-5mm \hspace{5mm}

How can we understand intelligent behavior? How can we design intelligent computers? These are questions that have been discussed by scientists and the public at large for over 50 years. As mathematicians, however, the question we want to ask is ``is there a {\it mathematical} theory underlying intelligence?'' I believe the first mathematical attack on these issues was Control Theory, led by Wiener and Pontryagin. They were studying how to design a controller which drives a motor affecting the world and also sits in a feedback loop receiving measurements from the world about the effect of the motor action. The goal was to control the motor so that the world, as measured, did something specific, i.e.\ move the tiller so that the boat stays on course. The main complication is that nothing is precisely predictable: the motor control is not exact, the world does unexpected things because of its complexities and the measurements you take of it are imprecise. All this led, in the simplest case, to a beautiful analysis known as the Wiener-Kalman-Bucy filter (to be described below).

But Control Theory is basically a theory of the output side of intelligence with the measurements modeled in the simplest possible way: e.g.\ linear functions of the state of the world system being controlled plus additive noise. The real input side of intelligence is perception in a much broader sense, the analysis of all the noisy incomplete signals which you can pick up from the world through natural or artificial senses. Such signals typically display a mix of distinctive patterns which tend to repeat with many kinds of variations and which are confused by noisy distortions and extraneous clutter. The interesting and important structure of the world is thus coded in these signals, using a code which is complex but not perversely so.

\subsection{Logic vs.\ Statistics}\vskip-5mm \hspace{5mm}

The first serious attack on problems of perception was the attempt
to recognize speech which was launched by the US defense agency
ARPA in 1970. At this point, there were two competing ideas of
what was the right formalism for combining the various clues and
features which the raw speech yielded. The first was to use logic
or, more precisely, a set of `production rules' to augment a
growing database of true propositions about the situation at hand.
This was often organized in a `blackboard', a two-dimensional
buffer with the time of the asserted proposition plotted along the
$x$-axis and the level of abstraction (i.e.\ signal --- phone ---
phoneme --- syllable --- word --- sentence) along the $y$-axis.
The second was to use statistics, that is, to compute
probabilities and conditional probabilities of various possible
events (like the identity of the phoneme being pronounced at some
instant). These statistics were computed by what was called the
`forward-backward' algorithm, making 2 passes in time, before the
final verdict about the most probable translation of the speech
into words was found. This issue of logic vs.\ statistics in the
modeling of thought has a long history going back to Aristotle
about which I have written in [M].

I think it is fair to say that statistics won. People in speech
were convinced in the 1970's, artificial intelligence researchers
converted during the 1980's as expert systems needed statistics so
clearly (see Pearl's influential book [P]), but vision researchers
were not converted until the 1990's when computers became powerful
enough to handle the much larger datasets and algorithms needed
for dealing with 2D images.

The biggest reason why it is hard to accept that statistics
underlies all our mental processes --- perception, thinking and
acting --- is that we are not consciously aware of 99\% of the
ambiguities with which we deal every second. What philosophers
call the `raw qualia', the actual sensations received, do not make
it to consciousness; what we are conscious of is a precise
unambiguous enhancement of the sensory signal in which our
expectations and our memories have been drawn upon to label and
complete each element of the percept. A very good example of this
comes from the psychophysical experiments of Warren \& Warren [W]
in 1970: they modified recorded speech by replacing a single
phoneme in a sentence by a noise and played this to subjects.
Remarkably, the subjects did {\it not} perceive that a phoneme was
missing but believed they had heard the one phoneme which made the
sentence semantically consistent:
\begin{center}
\begin{tabular} {l|l}
{\sc Actual sound} & {\sc Perceived words}\\
\hline
the ?`eel is on the shoe & the {\it h}eel is on the shoe\\
the ?`eel is on the car & the {\it wh}eel is on the car\\
the ?`eel is on the table & the {\it m}eal is on the table\\
the ?`eel is on the orange & the {\it p}eel is on the orange
\end{tabular}
\end{center}

\begin{figure}[htb]
\begin{center}
\epsfig{width=2.25in,file=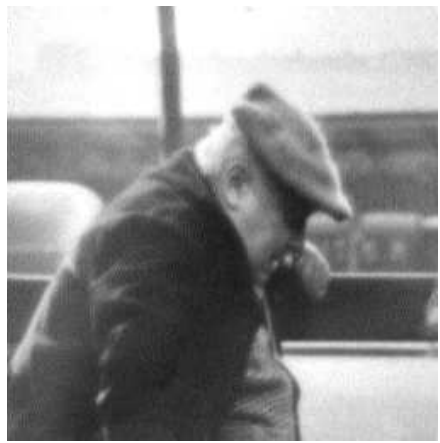}
\begin{minipage}[h]{10cm}
\caption{\footnotesize Why is this old man recognizable from a
cursory glance? His outline threads a complex path amongst the
cluttered background and is broken up by alternating highlights
and shadows and by the wrinkles on his coat. There is no single
part of this image which suggests a person unambiguously (the ear
comes closest but the rest of his face can only be guessed at). No
other object in the image stands out --- the man's cap, for
instance, could be virtually anything. Statistical methods, first
grouping contours, secondly guessing at likely illumination
effects and finally using probable models of clothes may draw him
out. No known computer algorithm comes close to finding a man in
this image.} \label{fig:oldman}
\end{minipage}
\end{center}
\end{figure}

Two things should be noted. Firstly, this showed clearly that the actual auditory signal did not reach consciousness. Secondly, the choice of percept was a matter of probability, not certainty. That is, one might find some odd shoe with a wheel on it, a car with a meal on it, a table with a peel on it, etc.\ but the words which popped into consciousness were the most likely. An example from vision of a simple image, whose contents require major statistical reasoning to reconstruct, is shown in figure~\ref{fig:oldman}.

It is important to clarify the role of probability in this approach. The uncertainty in a given situation need not be caused by observations of the world being truly unpredictable as in quantum mechanics or even effectively so as in chaotic phenomena. It is rather a matter of efficiency: in order to understand a sentence being spoken, we do not need to know all the things which affect the sound such as the exact acoustics of the room in which we are listening, nor are we even able to know other factors like the state of mind of the person we are listening to. In other words, we always have incomplete data about a situation. A vast number of physical and mental processes are going on around us, some germane to the meaning of the signal, some extraneous and just cluttering up the environment. In this `blooming, buzzing' world, as William James called it, we need to extract information and the best way to do it, apparently, is to make a stochastic model in which all the irrelevent events are given a simplified probability distribution. This is not unlike the stochastic approach to Navier-Stokes, where one seeks to replace turbulence or random molecular effects on small scales by stochastic perturbations.

\subsection{The Bayesian setup}\vskip-5mm \hspace{5mm}

Having accepted that we need to use probabilities to combine bits and pieces of evidence, what is the mathematical set up for this? We need the following ingredients: a) a set of random variables, some of which describe the observed signal and some the `hidden' variables describing the events and objects in the world which are causing this signal, b) a class of stochastic models which allow one to express the variability of the world and the noise present in the signals and c) specific parameters for the one stochastic model in this class which best describes the class of signals we are trying to decode now. More formally, we shall assume we have a set ${\bf x} = ({\bf x}_{\text{o}}, {\bf x}_{\text{h}})$ of observed and hidden random variables, which may have real values or discrete values in some finite or countable sets, we have a set ${\boldsymbol{\theta}}$ of parameters and we have a class of probability models $\text{Pr}({\bf x} \mid {\boldsymbol{\theta}})$ on the $x$'s for each set of values of the $\boldsymbol{\theta}$'s. The crafting or learning of this model may be called the first problem in the mathematical theory of perception. It is usual to factor these probability distributions:
$$ \text{Pr}({\bf x} \mid {\boldsymbol{\theta}}) = \text{Pr}({\bf x}_{\text{o}} \mid {\bf x}_{\text{h}}, \, \boldsymbol{\theta}) \cdot \text{Pr}({\bf x}_{\text{h}} \mid \boldsymbol{\theta}),$$
where the first factor, describing the likelihood of the observations from the hidden variables, is called the {\it imaging model} and the second, giving probabilities on the hidden variables, is called the {\it prior}. In the full Bayesian setting, one has an even stronger prior, a full probability model $\text{Pr}({\bf x}_h, \boldsymbol{\theta})$, including the parameters.

The second problem of perception is that we need to estimate the
values of the parameters $\boldsymbol{\theta}$ which give the best
stochastic model of this aspect of the  world. This often means
that you have some set of measurements $\{{\bf x}^{(\alpha)}\}$
and seek the value of $\boldsymbol{\theta}$ which maximizes their
likelihood $\prod_\alpha \text{Pr}({\bf x}^{(\alpha)} \mid
\boldsymbol{\theta)}$. If the hidden variables as well as the
observations are known, this is called supervised learning; if the
hidden variables are not known, then it is unsupervised and one
may maximize, for instance, $\prod_\alpha \sum_{{\bf x}_{\text h}}
\text{Pr}({\bf x}^{(\alpha)}_{\text o}, {\bf x}_{\text h} \mid
\boldsymbol{\theta})$. If one has a prior on the
$\boldsymbol{\theta}$'s too, one can also estimate them from the
mean or mode of the full posterior $\text{Pr}(\boldsymbol{\theta}
\mid \{ {\bf x}^{(\alpha)} \})$.

Usually a more challenging problem is how many parameters
$\boldsymbol{\theta}$ to include. At one extreme, there are simple
`off-the-shelf' models with very few parameters and, at the other
extreme, there are fully non-parametric models with infinitely
many parameters. Here the central issue is how much data one has:
for any set of data, models with too few parameters distort the
information the data contains and models with too many overfit the
accidents of this data set. This is called the {\it bias-variance
dilemma}. There are two main approaches to this issue. One is
cross-validation: hold back parts of the data, train the model to
have maximal likelihood on the training set and test it by
checking the likelihood of the held out data. There is also a
beautiful theoretical analysis of the problem due principally to
Vapnik [V] and involving the {\it VC dimension} of the models ---
the size of the largest set of data which can be split in all
possible ways into more and less likely parts by different choices
of $\boldsymbol{\theta}$.

As Grenander has emphasized, a very useful test for a class of models is to synthesize from it, i.e.\ choose random samples according to this probability measure and to see how well they resemble the signals we are accustomed to observing in the world. This is a stringent test as signals in the world usually express layers and layers of structure and the model tries to describe only a few of these.

The third problem of perception is using this machinary to actually perceive: we assume we have measured specific values ${\bf x}_{\text o} = \widehat{\bf x}_{\text o}$ and want to infer the values of the hidden variables ${\bf x}_{\text h}$ in this situation. Given these observations, by Bayes' rule, the hidden variables are distributed by the so-called {\it posterior} distribution:
$$ \text{Pr}({\bf x}_{\text h} \mid \widehat{\bf x}_{\text o},\, \boldsymbol{\theta}) =
\frac{\text{Pr}(\widehat{\bf x}_{\text{o}} \mid {\bf x}_{\text{h}},\, \boldsymbol{\theta}) \cdot \text{Pr}({\bf x}_{\text{h}} \mid \boldsymbol{\theta})}{\text{Pr}(\widehat{\bf x}_{\text o} \mid \boldsymbol{\theta})}\,
\propto \, \text{Pr}(\widehat{\bf x}_{\text{o}} \mid {\bf x}_{\text{h}},\, \boldsymbol{\theta}) \cdot \text{Pr}({\bf x}_{\text{h}} \mid \boldsymbol{\theta})
$$
One may then want to estimate the mode of the posterior, the most likely value of ${\bf x}_{\text h}$. Or one may want to estimate the mean of some functions $f({\bf x}_{\text h})$ of the hidden variables. Or, if the posterior is often multi-modal and some evidence is expected to available later, one usually wants a more complete description or approximation to the full posterior distribution.

\section{A basic example: HMM's and speech recognition}\label{section 2}
\setzero\vskip-5mm \hspace{5mm }

A convenient way to introduce the ideas of Pattern Theory is to outline the simple Hidden Markov Model method in speech recognition to illustrate many of the ideas and problems which occur almost everywhere. Here the observed random variables are  the values of the sound signal $s(t)$, a pressure wave in air. The hidden random variables are the states of the speaker's mouth and throat and the identity of the phonemes being spoken at each instant. Usually this is simplified, replacing the signal by samples $s_k = s(k \Delta t)$ and taking for hidden variables a sequence $x_k$ whose values indicate which phone in which phoneme is being pronounced at time $k \Delta t$. The stochastic model used is:
$$ \text{Pr}(x_\centerdot, s_\centerdot) = \prod_k p_1(x_k \mid x_{k-1}) p_2(s_k \mid x_k)$$
i.e.\ the $\{x_k\}$ form a Markov chain and each $s_k$ depends only on $x_k$. This is expressed by the graph:
$$ \begin{array}{ccccccc}
& \overset{s_{k-1}}{{\circ}} & & \overset{s_k}{{\circ}} & & \overset{s_{k+1}}{{\circ}} &\\
&\big\arrowvert&&\big\arrowvert&&\big\arrowvert&\\
\cdots\frac{\qquad}{\qquad}\hspace*{-.16in}&\underset{x_{k-1}}{\circ}&\hspace*{-.2in}\frac{\qquad}{\qquad}\hspace*{-.11in} & \underset{x_k}{\circ} & \hspace*{-.15in}\frac{\qquad}{\qquad}\hspace*{-.2in} & \underset{x_{k+1}}{\circ} & \hspace*{-.21in}\frac{\qquad}{\qquad}\cdots
\end{array} $$
in which each variable corresponds to a vertex and the graphical Markov property holds: if 2 vertices $a,b$ in the graph are separated by a subset $S$ of vertices, then the variables associated to $a$ and $b$ are conditionally independent if we fix the variables associated to $S$.

This simple model works moderately well to decode speech because of the linear nature of the graph, which allows the ideas of dynamic programming to be used to solve for the marginal distributions and the modes of the hidden variables, given any observations $\widehat{s}_\centerdot$. This is expressed simply in the recursive formulas:
\begin{eqnarray*}
\text{Pr}(x_k \mid \widehat{s}_{\le k}) &=& \frac{\sum_{x_{k-1}}
p_1(x_k \mid x_{k-1})p_2(\widehat{s}_k \mid x_k)
\text{Pr}(x_{k-1} \mid \widehat{s}_{\le (k-1)})}
{\sum_{x_k} \text{numerator}}\\
\text{maxp}(x_k , \widehat{s}_{\le k}) &\underset{\text{\scriptsize def}}{=}& \max_{x_{\le (k-1)}}
\text{Pr}(x_k, x_{\le k-1}, \widehat{s}_{\le k}) \\
&=& \max_{x_{k-1}} \left( p_1(x_k \mid x_{k-1}) p_2(\widehat{s}_x \mid x_k)
\text{maxp}(x_{k-1}, \widehat{s}_{\le (k-1)})\right).
\end{eqnarray*}
Note that if each $x_k$ can take $N$ values, the complexity of
each time step is $O(N^2)$.

In any model, if you can calculate the conditional probabilities of the hidden variables and if the model is of exponential type, i.e.\
$$ \text{Pr}(x_\centerdot \mid \boldsymbol{\theta}_\centerdot) =
\frac{1}{Z(\boldsymbol{\theta})}e^{\sum_k \theta_k \cdot E_k(x_\centerdot)},$$
then there is also an efficient method of optimizing the parameters $\boldsymbol{\theta}$. This is called the {\it EM algorithm} and, because it holds for HMM's, it is one of the key reasons for the early successes of the stochastic approach to speech recognition. For instance, a Markov chain $\{x_k\}$ is an exponential model if we let the $\boldsymbol{\theta}$'s be $\log(p(a \mid b))$ and write the chain probabilities as:
$$ \text{Pr}(x_\centerdot) = e^{\sum_{a,b} \log(p(a \mid b) |\{k \mid x_k=a, x_{k-1} = b\}|}.$$
The fundamental result on exponential models is that the $\boldsymbol{\theta}$'s are determined by the expectations $\widehat{E}_k=\text{Exp}(E_k)$ and that any set of expectations $\widehat{E}_k$ that can be achieved in some probability model (with all probabilities non-zero), is also achieved in an exponential model.

\subsection{Continuous and discrete variables}\vskip-5mm \hspace{5mm}

In this model, the observations $s_k$ are naturally continuous random variables, like all primary measurements of the physical world. But the hidden variables are discrete: the set of phonemes, although somewhat variable from language to language, is always a small discrete set. This combination of discrete and continuous is characteristic of perception. It is certainly a psychophysical reality: for example experiments show that our perceptions lock onto one or another phoneme, resisting ambiguity (see [L], Ch.8, esp.\ p.176). But it shows itself more objectively in the low-level statistics of natural signals. Take almost any class of continuous real-valued signals $s(t)$ generated by the world and compile a histogram of their changes $x=s(t+\Delta t) - s(t)$ over some fixed time interval $\Delta t$. This empirical distribution will very likely have kurtosis ($= \text{Exp}\left((x-\bar{x})^4\right)/\sigma(x)^4$) greater than 3, the kurtosis of any Gaussian distribution! This means that, compared to a Gaussian distribution with the same mean and standard deviation, $x$ has higher probability of being quite small or quite large but a lower probability of being average. Thus, compared to Brownian motion, $s(t)$ tends to move relatively little most of the time but to make quite large moves sometimes. This can be made precise by the theory of stochastic processes with iid increments, a natural first approximation to any stationary Markov process. The theory of such processes says that (a) their increments always have kurtosis at least 3, (b) if it equals 3 the process is Brownian and (c) if it is greater, samples from the process almost surely have discontinuities. At the risk of over-simplfying, we can say {\it kurtosis $>$ 3 is nature's universal signal of the presence of discrete events/objects in continuous space-time}.

A classic example of this are stock market prices. Their changes (or better, changes in log(price)) have a highly non-Gaussian distribution with polynomial tails. In speech, the changes in the log(power) of the windowed Fourier transform show the same phenomenon, confirming that $s(t)$ cannot be decently modeled by colored Gaussian noise.

\subsection{When compiling full probability tables is impractical}\vskip-5mm \hspace{5mm}

Applying HMM's in realistic settings, it usually happens that $N$
is too large for an exhaustive search of complexity $O(N^2)$ or
that the $x_k$ are real valued and, when adequately sampled, again
$N$ is too large. There is one other situation in which the
HMM-style approach works easily --- the {\it Kalman filter}. In
Kalman's setting, each variable $x_k$ and $s_k$ is real
vector-valued instead of being discrete and $p_1$ and $p_2$ are
{\it Gaussian distributions with fixed covariances and means
depending linearly on the conditioning variable}. It is then easy
to derive recursive update formulas, similar to those above, for
the conditional distributions on each $x_k$, given the past data
$\widehat{s}_{\le k}$.

But usually, in the real-valued variable setting, the $p$'s are more complex than Gaussian distributions. An example is the tracking problem in vision: the position and velocity $x_k$ of some specific moving object at time $k\Delta t$ is to be inferred from a movie $\widehat{s}_k$, in which the object's location is confused by clutter and noise. It is clear that the search for the optimal reconstruction $x_k$ must be pruned or approximated. A dramatic breakthrough in this and other complex situations has been to adapt the HMM/Kalman ideas by using weak approximations to the marginals $\text{Pr}(x_k \mid \widehat{s}_{\le k})$ by a finite set of samples, an idea called {\it particle filtering}:
\begin{eqnarray*}
\text{Pr}(x_k \mid \widehat{s}_{\le k}) &\underset{\text{\footnotesize weak}}{\sim}& \sum_{i=1}^N w_{i,k} \delta_{x_{i,k}} (x_k),
\quad \text{that is,}\\
\text{Exp}(f(x_k) \mid \widehat{s}_{\le k}) &\approx& \sum_{i=1}^N w_{i,k} f(x_{i,k}), \text{ for suitable } f.
\end{eqnarray*}
This idea was proposed originally by Gordon, Salmond and Smith [G-S-S] and is developed at length in the recent survey [D-F-G]. An example with explicit estimates of the posterior from the work of Isard and Blake [I-B] is shown in figure 2. They follow the version known as bootstrap particle filtering in which, for each $k$, $N$ samples $x'_l$ are drawn with replacement from the weak approximation above, each sample is propagated randomly to a new sample $x''_l$ at time $(k+1)$ using the prior $p(x_{k+1} \mid x'_l)$ and these are reweighted proportional to $p(\widehat{s}_{k+1} \mid x''_l)$.
\begin{figure}[htb]
%\centerline{\epsfig{width=3.5in,file=Blake.eps}}
\begin{center}
\epsfig{width=3.5in,file=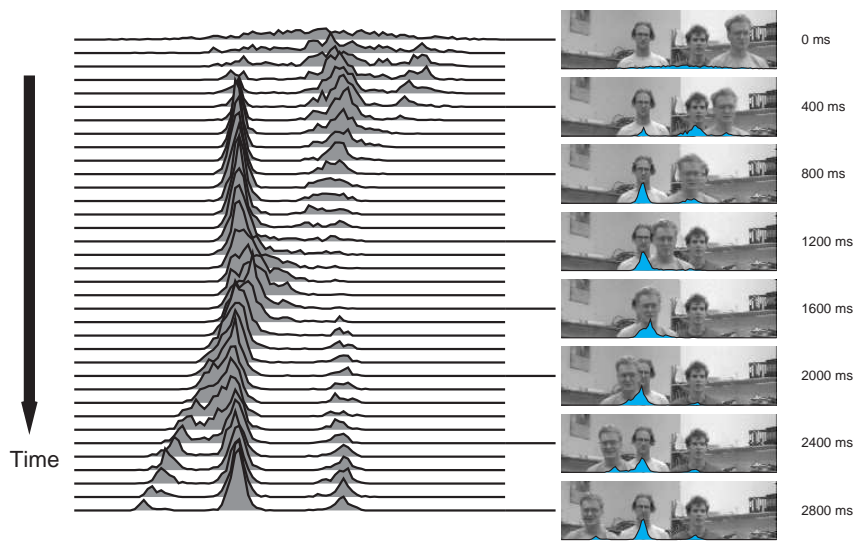}
\begin{minipage}[h]{10cm}
 \caption{\footnotesize Work
of Blake and Isard tracking three faces in a moving image
sequence. The curves represent estimates of the posterior
probability distributions for faces at each location obtained by
smoothing the weighted sum of delta functions at the `particles'.
Note how multi-modal these are and how the tracker recovers from
the temporary occlusion of one face by another.}
\end{minipage}
\end{center}
\end{figure}

\subsection{No process in nature is truly Markov}\vskip-5mm \hspace{5mm}

A more serious problem with the HMM approach is that the Markov assumption is never really valid and it may be much too crude an approximation. Consider speech recognition. The finite lexicon of words clearly constrains the expected phoneme sequences, i.e.\ if $x_k$ are the phonemes, then $p_1(x_k \mid x_{k-1})$ depends on the current word(s) containing these phonemes, i.e.\ on a short but variable part of the preceding string $\{x_{k-1}, x_{k-2}, \cdots \}$ of phonemes. To fix this, we could let $x_k$ be a pair consisting of a word and a specific phoneme in this word; then $p_1(x_k \mid x_{k-1})$ would have two quite different values depending on whether $x_{k-1}$ was the last phoneme in the word or not. Within a word, the chain needs only to take into account the variability with which the word can be pronounced. At word boundaries, it should use the conditional probabilities of word pairs. This builds much more of the patterns of the language into the model.

Why stop here? State-of-the-art speech recognizers go further and let $x_k$ be a pair of consecutive words plus a triphone\footnote{So-called co-articulation effects mean that the pronunciation of a phoneme is affected by the preceding and suceeding phonemes.} in the second word (or bridging the first and second word) whose middle phoneme is being pronounced at time $k\Delta t$. Then the transition probabilities in the HMM involve the statistics of `trigrams', consecutive word triples in the language. But grammar tells us that words sequences are also structured into phrases and clauses of variable length forming a parse tree. These clearly affect the statistics. Semantics tells us that words sequences are further constrained by semantic plausibility (`sky' is more probable as the word following `blue' than `cry') and pragmatics tells us that sentences are part of human communications which further constrain probable word sequences.

All these effects make it clear that certain parts of the signal should be grouped together into units on a higher level and given labels which determine how likely they are to follow each other or combine in any way. This is the essence of grammar: higher order random variables are needed whose values are subsets of the low order random variables. The simplest class of stochastic models which incorporate variable length random substrings of the phoneme sequence are {\it probabilistic context free grammars} or PCFG's. Mathematically, they are a particular type of random branching tree.

{\bf Definition} {\it A \underline{PCFG} is a stochastic model in which the random variables are (a) a sequence of rooted trees $\{{\cal T}_n\}$, (b) a linearly ordered sequence of observations $s_k$ and a 1:1 correspondence between the observations $s_k$ and the leaves of the whole forest of trees such that the children of any vertex of any tree form an interval $\{s_{k}, s_{k+1}, \cdots, s_{k'}\}$ in time and (c) a set of labels $x_v$ for each vertex. The probability model is given by conditional probabilities $p_1(x_{v_k} \mid x_v)$ for the labels of each child of each vertex\footnote{Caution to specialists: our label $x_v$ is the name of the `production rule' with this vertex as its head, esp.\ it fixes the arity of the vertex. We are doing it this way to simplify the Markov property.} and $p_2(s_k \mid x_{v_k})$ for the observations, conditional on the label of the corresponding leaf.}

See figure 3 for an example. This has a Markov property if we define the `extended' state $x_k^*$ at leaf $k$ to be not only the label $x_k$ at this leaf but the whole sequence of labels on the path from this leaf to the root of the tree in which this leaf lies. Conditional on this state, the past and the future are independent.

This is a mathematically elegant and satisfying theory: unfortunately, it also fails, or rather explodes because, in carrying it out, the set of labels gets bigger and bigger. For instance, it is not enough to have a label for noun phrase which expands into an adjective plus a noun. The adjective and noun must agree in number and (in many languages) gender, a constraint that must be carried from the adjective to the noun (which need not be adjacent) via the label of the parent. So we need 4 labels, all combinations of singular/plural masculine/feminine noun phrases. And semantic constraints, such as Pr(`blue sky') $>$ Pr(`blue cry'), would seem to require even more labels like `colorable noun phrases'. Rather than letting the label set explode, it is better to consider a bigger class of grammars, which express these relations more succinctly but which are not so easily converted into HMM's: {\it unification grammars} [Sh] or {\it compositional grammars} [B-G-P]. The need for grammars of this type is especially clear when we look at formalisms for expressing the grouping laws in vision: see figure 3. The further development of stochastic compositional grammars, both in language and vision, is one of the main challenges today.
\begin{center}
\epsfig{width=4.5in,file=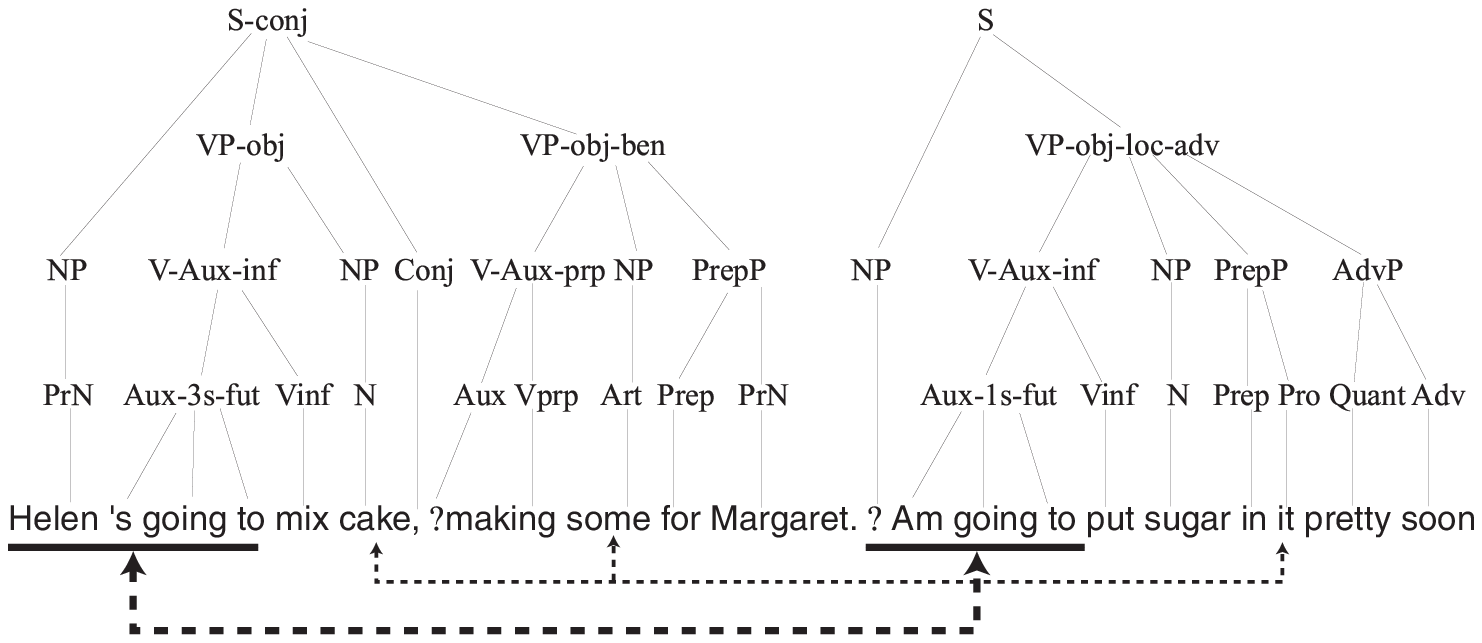} \\
\vspace{.1in}
\epsfig{width=2.2in,file=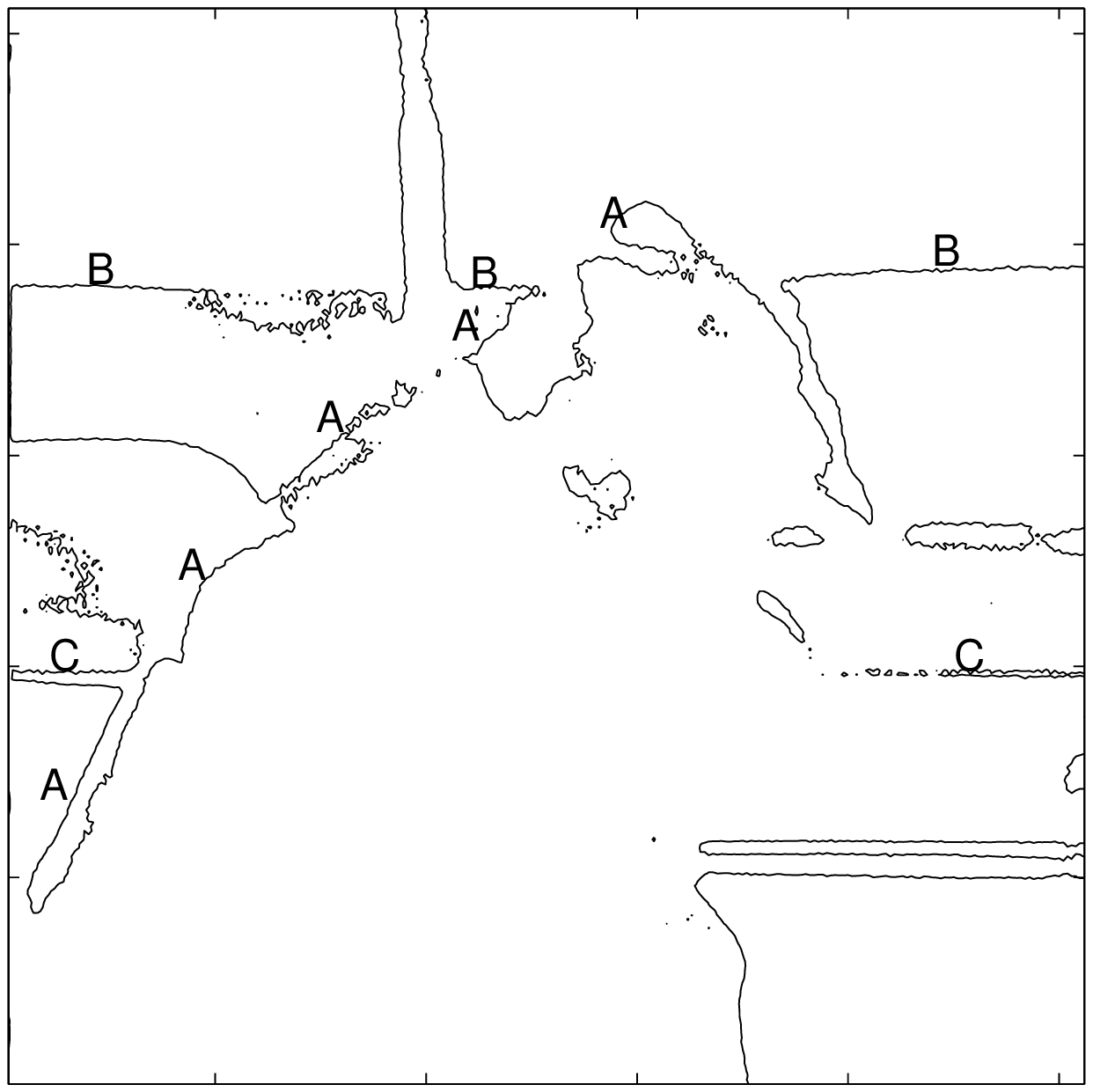} \vspace*{-.2in}
\end{center}
\begin{figure}[h]%[htb]
\begin{center}
\begin{minipage}{10cm}
\caption{\footnotesize Grouping in language and vision: On top,
parsing the not quite grammatical speech of a 2 1/2 year old Helen
describing her own intentions ([H]): above the sentence, a
context-free parse tree; below it, longer range non-Markov links
--- the identity `cake'=`some'=`it' and the unification of the two
parts `Helen's going to' = `(I) am going to'. On the bottom, 2
kinds of grouping with an iso-intensity contour of the image in
Figure 1: note the broken but visible contour of the back marked
by `A' and the occluded contours marked by `B' and `C' behind the
man.}
\end{minipage}
\end{center}
\end{figure}

\section{The `natural degree of generality': MRF's or Graphical Models}
\label{section 3}
\setzero\vskip-5mm \hspace{5mm }

The theory of HMM's deals with one-dimensional signals. But
images, the signals occurring in vision, are usually
two-dimensional --- or three-dimensional for MR scans and movies
(3 space dimensions and 2 space plus 1 time dimension), even
four-dimensional for echo cardiograms. On the other hand, the
parse tree is a more abstract graphical structure and other
`signals', like medical data gathered about a particular patient,
are structured in complex ways (e.g.\ a set of blood tests, a
medical history). This leads to the basic insight of Grenander's
Pattern Theory [G]: that the variables describing the structures
in the world are typically related in a graphical fashion, edges
connecting variables which have direct bearing on each other.
Finding the right graph or class of graphs is a crucial step in
setting up a satisfactory model for any type of patterns. Thus the
applications, as well as the mathematical desire to find the most
general setting for this theory, lead to the idea of replacing a
simple chain of variables by a set of variables with a more
general graphical structure. The general concept we need is that
of a Markov random field:

{\bf Definition} {\it A \underline{Markov random field} is a graph $G = (V,E)$, a set of random variables $\{x_v\}_{v \in V}$, one for each vertex, and a joint probability distribution on these variables of the form:
$$ \mathrm{Pr}(x_\centerdot) = \frac{1}{Z} e^{-\sum_C E_C(\{x_v\}_{v \in C})},$$
where $C$ ranges over the cliques (fully connected subsets) of the graph, $E_C$ are any functions and $Z$ a constant. If the variables $x_v$ are real-valued for $v \in V'$, we make this into a probability density, multiplying by $\prod_{n \in V'} dx_n$. Moreover, we can put each model in a family by introducing a temperature $T$ and defining:}
$$ \mathrm{Pr}_T(x_\centerdot) = \frac{1}{Z_T} e^{-\sum_C E_C(\{x_v\}_{v \in C})/T}.$$

These are also called {\it Gibbs models} in statistical mechanics (where the $E_C$ are called {\it energies}) and {\it graphical models} in learning theory and, like Markov chains, are characterized by their conditional independence properties. This characterization, called the Hammersley-Clifford theorem, is that if two vertices $a,b \in V$ are separated by a subset $S \subset V$ (all paths in $G$ from $a$ to $b$ must include some vertex in $S$), then $x_a$ and $x_b$ are conditionally independent given $\{x_v\}_{v \in S}$. The equivalence of these independence properties, plus the requirement that all probabilities be positive, with the simple explicit formula for the joint probabilities makes it very convincing that MRF's are a natural class of stochastic models.

\subsection{The Ising model}\vskip-5mm \hspace{5mm}

This class of models is very expressive and many types of patterns which occur in the signals of nature can be captured by this sort of stochastic model. A basic example is the Ising model and its application to the image segmentation problem. In the simplest form, we take the graph $G$ to be a square $N \times N$ grid with two layers, with observable random variables $p_{i,j} \in \mathbb{R}, 1 \le i,j \le N$ associated to the top layer and hidden random variables $x_{i,j} \in \{+1, -1\}$ associated to the bottom layer. We connect by edges each $x_{i,j}$ vertex to the $p_{i,j}$ vertex above it and to its 4 neighbors $x_{i \pm 1,j}, x_{i,j \pm 1}$ in the $x$-grid (except when the neighbor is off the grid) and no others. The cliques are just the pairs of vertices connected by edges. Finally, we take for energies:
\begin{eqnarray*}
E_C &=& -x_{i,j} \cdot x_{i',j'}, \text{ when } C = \{(i,j),(i',j')\}, \text{ two adjacent vertices in the $x$-grid},\\
E_C &=& -x_{i,j} \cdot y_{i,j}, \text{ when $C$ consists of the $(i,j)$ vertices in the $x$- and $y$-grids}.
\end{eqnarray*}
The modes of the posteriors $\text{Pr}_T(x_\centerdot \mid \widehat{y}_\centerdot)$ are quite subtle: $x$'s at adjacent vertices try to be equal but they also seek to have the same sign as the correponding $\widehat{y}$. If $\widehat{y}$ has rapid positive and negative swings, these are in conflict. Hence the more probable values of $x$ will align with the larger areas where $\widehat{y}$ is consistently of one sign. This can be used to model a basic problem in vision: the {\it segmentation problem}. The vision problem is to decompose the domain of an image $y$ into parts where distinct objects are seen. For example, the oldman image might be decomposed into 6 parts: his body, his head, his cap, the bench, the wall behind him and the sky. The decomposition is to be based on the idea that the image will tend to either slowly varying or to be statistically stationary at points on one object, but to change abruptly at the edges of objects. As proposed in [G-G], the Ising model can be used to treat the case where the image has 2 parts, one lighter and one darker, so that at the mode of the posterior the hidden variables $x$ will be $+1$ on one part, $-1$ on the other. An example is shown in figure 4. This approach makes a beautiful link between statistical mechanics and perception, in which the process of finding global patterns in a signal is like forming large scale structures in a physical material as the temperature cools through a phase transition.

\begin{figure}[hbp]
\begin{center}
\epsfig{width=1.65in,file=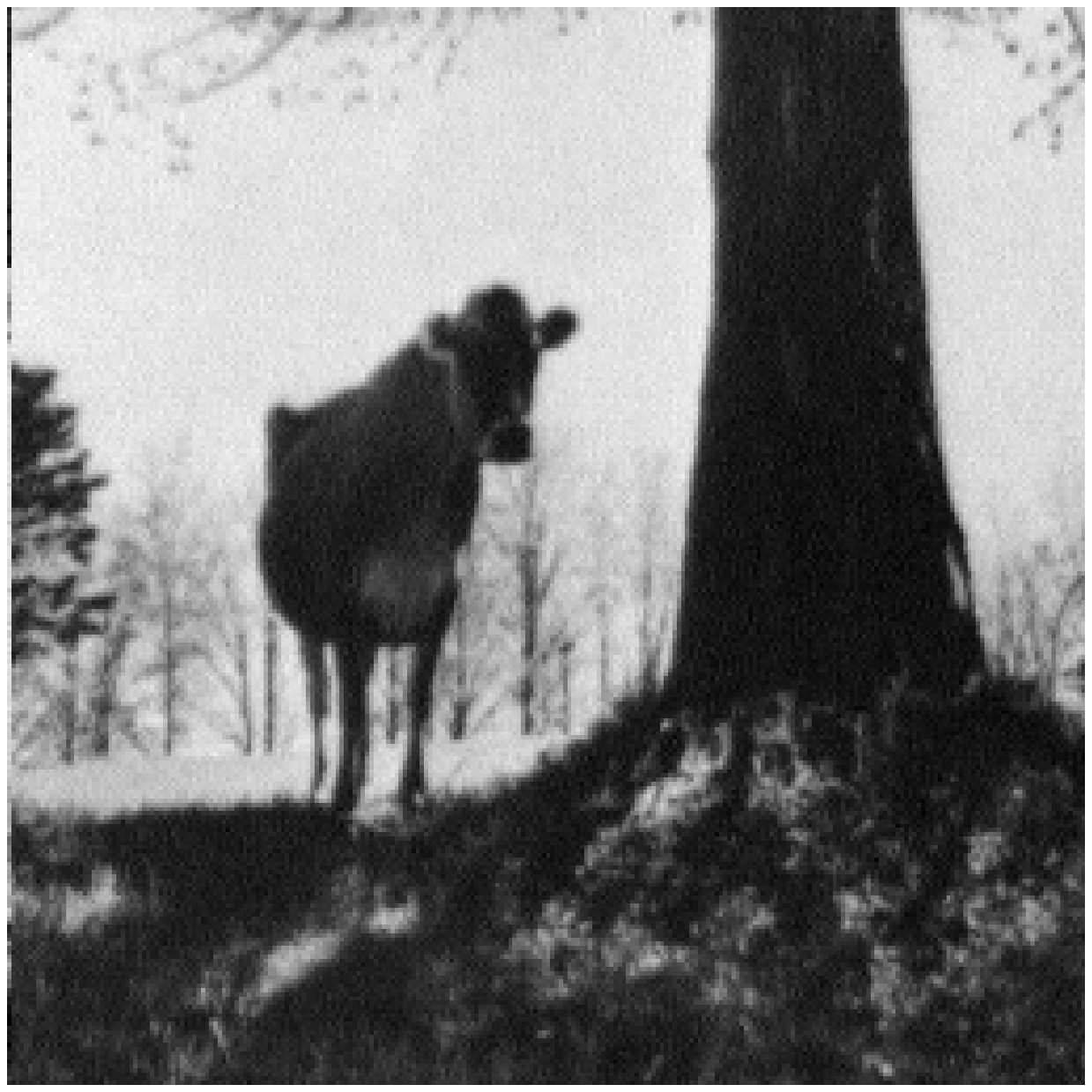} \hspace*{.5in} \epsfig{width=1.65in,file=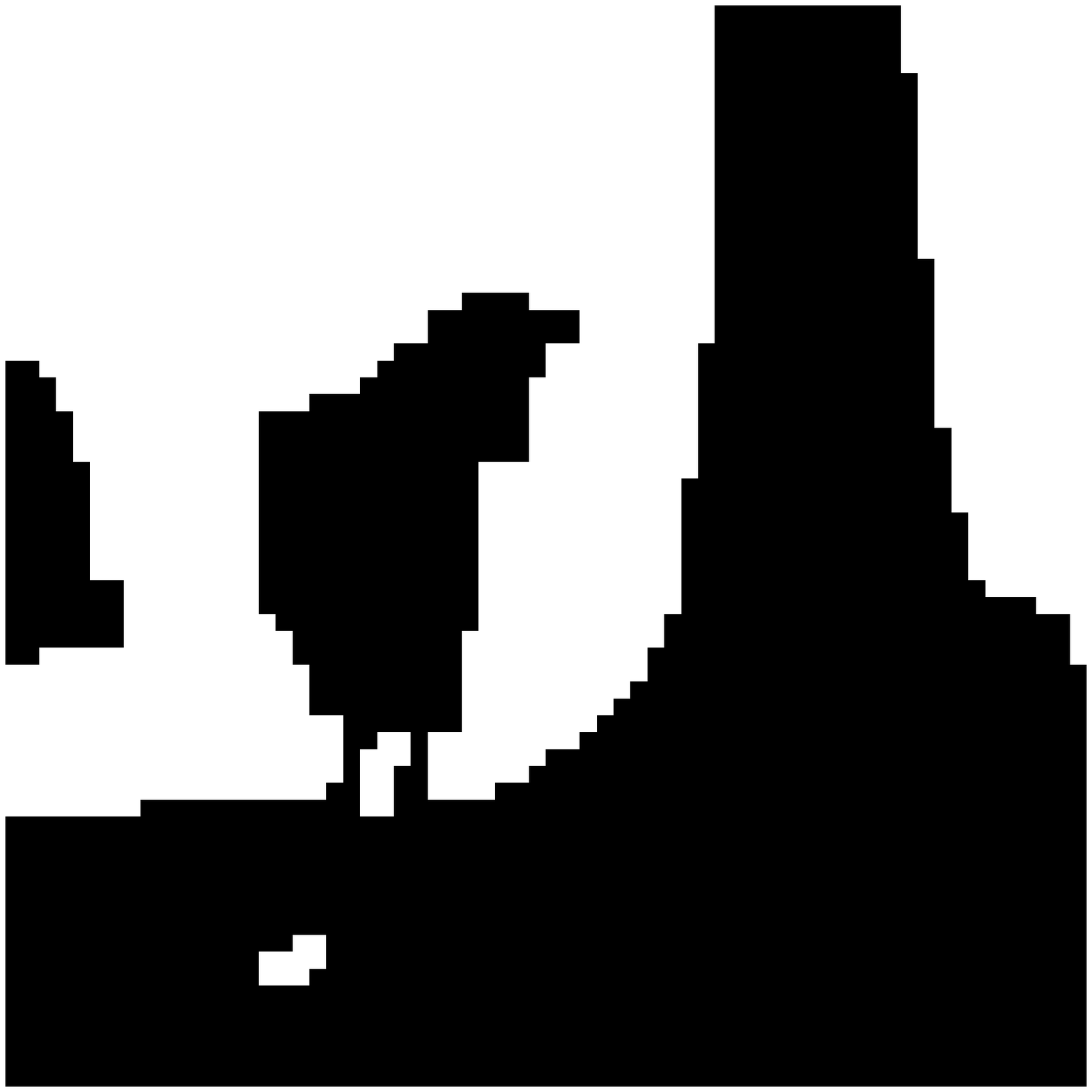}
\vspace*{.25in}\\
\epsfig{width=5in,file=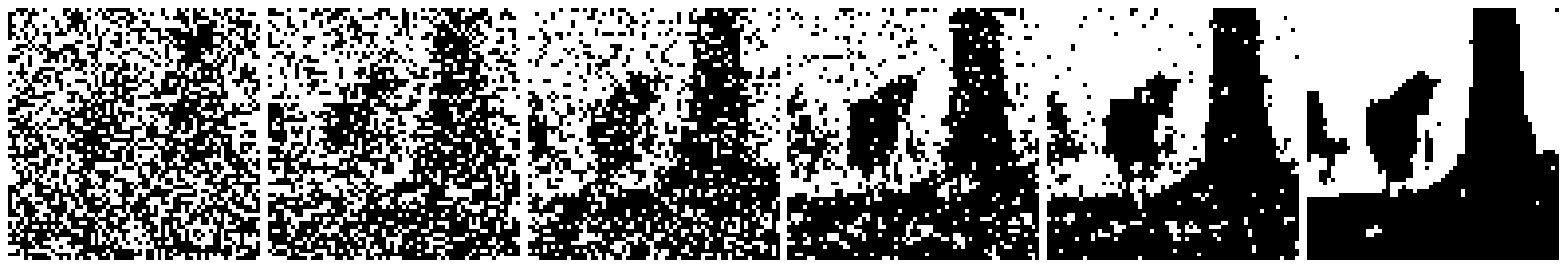}
\begin{minipage}[h]{10cm}
\caption{\footnotesize Statistical mechanics can be applied to the segmentation of images. On the top left, a rural scene taken as the external magnetic field, with its intensity scaled so that dark areas are negative, light areas are positive. At the top right, the mode or ground state of the Ising model. Along the bottom, the Gibbs distribution is sampled at a decreasing sequence of temperatures, discovering the global pattern bit by bit.}
\end{minipage}
\end{center}
\end{figure}

More complex models of this sort have been used extensively in image analysis, for texture segmentation, for finding disparity in stereo vision, for finding optic flow in moving images and for finding other kinds of groupings. We want to give one example of the expressivity of these models which is quite instructive. We saw above that exponential models can be crafted to reproduce some set of observed expectations but we also saw that scalar statistics from natural signals typically have high kurtosis, i.e.\ significant outliers, so that their whole distribution and not just their mean needs to be captured in the model. Putting these 2 facts together suggests that we seek exponential models which duplicate the whole distribution of some important statistics $f_\centerdot$.  This can be done using as parameters not just unknown constants but unknown functions:
$$ \text{Pr}(x_\centerdot \mid \phi_\centerdot) = \frac{1}{Z(\vec{\phi_\centerdot})}
e^{\sum_k \phi_k(f_k(x_\centerdot))}.$$
If $f_k$ depends only the variables $x_v \in C_k$, for some clique $C_k$, this is a MRF, whose energies have unknown functions in them. An example of this fitting is shown in Figure 5.

\begin{figure}[htp]
%\centerline{\epsfig{width=4in,file=TextEx.eps}}
\begin{center}
\epsfig{width=4in,file=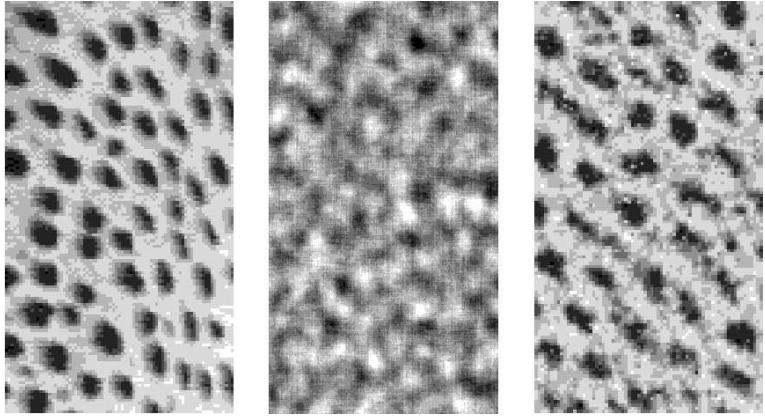}
\begin{minipage}[h]{10cm}
\caption{\footnotesize On the left, an image of the texture of a
Cheetah's hide, in the middle a synthetic image from the Gaussian
model with the same second order statistics, on the right a
synthetic image in which the full distribution on 7 filter
statistics are reproduced by an exponential model.}
\end{minipage}
\end{center}
\end{figure}

\subsection{Bayesian belief propagation}\vskip-5mm \hspace{5mm}

However, a problem with MRF models is that the dynamic programming
style algorithm used in speech and one-dimensional models to find
the posterior mode has no analog in 2D. One strategy for dealing
with this, which goes back to Metropolis, is to imitate physics
and introduce an artifical dynamics into the state space whose
equilibrium is the Gibbs distribution. This dynamics is called a
{\it Monte Carlo Markov Chain} (MCMC) and is how the panels in
figure 4 were generated. Letting the temperature converge to zero,
we get {\it simulated annealing} (see [G-G]) and, if we do it
slowly enough, will find the mode of the MRF model. Although slow,
this can be speeded up by biasing the dynamics (called {\it
importance sampling} --- see [T-Z] for a state-of-the-art
implementation with many improvements) and is an important tool.

Recently, however, another idea due to Weiss and collaborators (see [Y-F-W]) and linked to statistical mechanics has been found to give new and remarkably effective algorithms for finding these modes. From an algorithmic standpoint, the idea is to use the natural generalization of dynamic programming, called {\it Bayesian Belief Propagation} (BBP), which computes the marginals and modes correctly whenever the graph is a tree and just use it anyway on an arbitrary graph $G$! Mathematically, it amounts to working on the universal covering graph $\widetilde{G}$, which is a tree, hence much simpler, instead of $G$. In statistical mechanics, this idea is called the {\it Bethe approximation}, introduced by him in the 30's.

To explain the idea, start with the {\it mean field approximation}. The mean field idea is to find the best approximation of the MRF $p$ by a probability distribution in which the variables $x_v$ are all independent. This is formulated as the distribution $\prod_v p_v(x_v)$ which minimizes the Kullback-Liebler divergence $\text{KL}(\prod_v p_v, p)$. Unlike computing the true marginals of $p$ on each $x_v$ which is very hard, this approximation can be found by solving iteratively a coupled set of non-linear equations for the $p_v$. But the assumption of independence is much too restrictive. The idea of Bethe is instead to approximate $p$ by a $\pi_1(G)$-invariant distribution on $\widetilde{G}$.

Such distributions are easy to describe: note that a Markov random field on a {\it tree} is uniquely determined by its marginals $p_e(x_v,x_w)$ for each edge $e =(v,w)$ and, conversely, if we are given a compatible set of distributions $p_e$ for each edge (in the sense that, for all edges $(v,w_k)$ abutting a vertex $v$, the marginals of $p_{(v,w_k)}$ give distributions on $v$ independent of $k$), they define an MRF on $G$. So if we start with a Markov random field on any $G$, we get a $\pi_1(G)$-invariant Markov random field on $\widetilde{G}$ by making duplicate copies for each random variable $x_v, v \in V$ for each $\widetilde{v} \in \widetilde{V}$ over $v$ and lifting the edge marginals. But more generally, if we have any compatible set of probability distributions $\{p_e(v,w)\}_{ e \in E}$ on $G$, we also get a $\pi_1(G)$-invariant MRF on $\widetilde{G}$. Then the Bethe approximation is that family $\{p_e\}$ which minimizes $\text{KL}(\{p_e\},p)$. As in the mean field case, there is a natural iterative method of solving for this minimum, which turns out, remarkably, to be identical to the generalization of BBP to general graphs $G$.

This approach has proved effective in some cases at finding best segmentations of images via the mode of a two-dimensional MRF. Other interesting ideas have been proposed for solving the segmentation problem which we do not have time to sketch:  region growing, see esp.\ [Z-Y]), using the eigenfunctions of the graph-theoretic Laplacian, see [S-M], and multi-scale algorithms, see [P-B] and [S-B-B].

\section{Continuous space and time and continuous sets of random variables}
\label{section 4}
\setzero\vskip-5mm \hspace{5mm }

Although signals as we measure them are always sampled discretely, in the world itself signals are functions on the continua, time or space or both together. In some situations, a much richer mathematical theory emerges by replacing a countable collection of random variables by random processes and asking whether we can find good stochastic models for these continuous signals. I want to conclude this talk by mentioning three instances where some interesting analysis has arisen when passing to the continuum limit and going into some detail on two. We will not worry about algorithmic issues for these models.

\subsection{Deblurring and denoising of images}\vskip-5mm \hspace{5mm}

This is the area where the most work has been done, both because of its links with other areas of analysis and because it is one of the central problems of image processing. You observe a degraded image $I(x,y)$ as a function of continuous variables and seek to restore it, removing simultaneously noise and blur. In the discrete setting, the Ising model or variants thereof discussed above can be applied for this. There are two closely related ways to pass to the continuous limit and reformulate this as a problem in analysis. As both drop the stochastic interpretation and have excellent treatments in the literature, we only mention briefly one of a family of variants of each approach:

{\it Optimal piecewise smooth approximation of I via a variational problem}:
$$\min_{J,\Gamma} \left(c_1 \iint_{\!\! D} (I-J)^2 dxdy + c_2 \iint_{\!\! D-\Gamma} \|\nabla J\|^2 dxdy +c_3 |\Gamma| \right)$$
where $J$, the improved image, has discontinuities along the set of `edge' curves $\Gamma$. This approach is due to the author and Shah and has been extensively pursued by the schools of DeGiorgi and Morel. See [M-S]. It is remarkable that it is still unknown whether the minima to this functional are well behaved, e.g.\ whether $\Gamma$ has a finite number of components. Stochastic variants of this approach should exist.

{\it Non-linear diffusion of I}:
$$\frac{\partial J}{\partial t} =
\text{div}\left(\frac{\nabla J}{\|\nabla J\|}\right)+ \lambda
(I-J)$$ where $J$ at some future time is the enhancement. This
approach started with the work of Perona and Malik and has been
extensively pursued by Osher and his coworkers. See [Gu-M]. It can
be interpreted as gradient descent for a variant of the previous
variational problem.

\pagebreak
\subsection{Self-similarity of image statistics and image models}\vskip-5mm \hspace{5mm}

One of the earliest discoveries about the statistics of images $I$ was that their power spectra tend to obey power laws $$\text{Exp}|(\widehat{I}(\xi,\eta)|^2 \approx
(\xi^2 + \eta^2)^{-\lambda/2},$$
where $\lambda$ varies somewhat from image to image but clusters around the value 2. This has a very provocative interpretation: this power law is implied by self-similarity! In the language of lattice field theory, if $I(i,j), i,j \in \mathbb{Z}$ is a random lattice field and $\bar{I}$ is the block averaged field
$$\bar{I}(i,j) = \frac{1}{4} \left(I(2i,2j)+I(2i+1,2j)+I(2i,2j+1)+I(2i+1,2j+1)\right),$$
then we say the field is a renormalization fixed point if the distributions of $I$ and of $\bar{I}$ are the same. The hypothesis that natural images of the world, treated as a single large database, have renormalization invariant statistics has received remarkable confirmation from many quite distinct tests.

Why does this hold? It certainly isn't true for auditory or tactile signals. I think there is one major and one minor reason for it. The major one is that the world is viewed from a random viewpoint, so one can move closer or farther from any scene. To first approximation, this scales the image (though not exactly because nearer objects scale faster than distant ones). The minor one is that most objects are opaque but have, by and large, parts or patterns on them and, in turn, belong to clusters of larger things. This observation may be formulated as saying the world is not merely made up of objects but it is cluttered with them.

The natural setting for scale invariance is pass to the limit and model images as random functions $I(x,y)$ of two real variables. Then the hypothesis is that a suitable function space supports a probability measure which is invariant under both translations and scalings $(x,y) \mapsto (\sigma x, \sigma y)$, whose samples are `natural images'. This hypothesis encounters, however, an infra-red and an ultra-violet catastrophe:\\
a) The infra-red one is caused by larger and larger scale effects giving bigger and bigger positive and negative swings to a local value of $I$. But these large scale effects are very low-frequency and this is solved by considering $I$ to be defined only modulo an unknown constant, i.e.\ it is a sample from a measure on a function space mod constants.\\
b) The ultra-violet one is worse: there are more and more local oscillations of the signals at finer and finer scales and this contradicts Lusin's theorem that an integrable function is continuous outside sets of arbitrarily small measure. In fact, it is a theorem that {\it there is no translation and scale invariant probability measure on the space of locally integrable functions mod constants}. This can be avoided by allowing images to be generalized functions. In fact, the support can be as small as the intersection of all negative Sobolev spaces $\bigcap_\epsilon {\cal H}^{-\epsilon}$.

To summarize what a good statistical theory of natural images should explain, we have scale-invariance as just described, kurtosis greater than 3 as described in section 2.1 and finally the right local properties:
\begin{description}
\item[{\bf Hypothesis I}] A theory of images is a translation and scale invariant probability measure on the space of generalized functions $I(x,y)$ mod constants.
\item[{\bf Hypothesis II}] For any filter $F$ with mean 0, the marginal statistics of $F \ast I(x,y)$ have kurtosis greater than 3.
\item[{\bf Hypothesis III}] The local statistics of images reflect the preferred local geometries, esp.\ images of straight edges, but also curved edges, corners, bars, `T-junctions' and `blobs' as well as images without geometry, blank `blue sky' patches.
\end{description}
Hypothesis III is roughly the existence of what Marr, thinking globally of the image called the {\it primal sketch} and what Julesz, thinking locally of the elements of texture, referred to as {\it textons}. By scale invariance, the local and global image should have the same elements.

To quantify Hypothesis III, what is needed is a major effort at data mining. Specifically, the natural approach seems to be to take a small filter bank of zero mean local filters $F_1,\cdots,F_k$, a large data base of natural images $I_\alpha$ leading to the sample of points in $\mathbb{R}^k$ given by $(F_1 \ast I_\alpha(x,y), \cdots, F_K \ast I_\alpha(x,y)) \in \mathbb{R}^k$ for all $\alpha,x \text{ and } y$. One seeks a good non-parametric fit to this dataset. But Hypothesis III shows that this distribution will not be simple. For example Lee et al [L-P-M] have taken $k=8$, $F_i$ a basis of zero mean filters with fixed $3 \times 3$ support. They then make a linear tranformation in $\mathbb{R}^8$ normalizing the covariance of the data to $I_8$ (`whitening' the data), and to investigate the outliers, map the data with norms in the upper 20\% to $S^7$ by dividing by the norm. The analysis reveals that the resulting data has asymptotic infinite density along a non-linear surface in $S^7$! This surface is constructed by starting with an ideal image, black and white on the two sides of a straight edge and forming a $3 \times 3$ discrete image patch by integrating this ideal image over a tic-tac-toe board of square pixels. As the angle of the edge and the offset of the pixels to the edge vary, the resulting patches form this surface. This is the most concrete piece of evidence showing the complexity of local image statistics.

Are there models for these three hypotheses? We can satisfy the
first hypothesis by the unique scale-invariant Gaussian model,
called the free field by physicists --- but its samples look like
clouds and its marginals have kurtosis 3, so neither the second
nor third hypothesis is satisfied. The next best approximation
seems to be to use infinitely divisible measures, such as the
model constructed by the author and B.Gidas [M-G], which we call
{\it random wavelet expansions}:
$$
I(x,y) = \sum_i \phi_i(e^{r_i}x - x_i, e^{r_i}y - y_i),
$$
where $\{(x_i,y_i,r_i)\}$ is a Poisson process in $\mathbb{R}^3$ and $\phi_i$ are samples from an auxiliary Levi measure, playing the role of individual random wavelet primitives. But this model is based on adding primitives, as in a world of transparent objects, which causes the probability density functions of its marginal filter statistics to be smooth at 0 instead of having peaks there, i.e.\ the model does not produce enough `blue sky' patches with very low constrast.

A better approach are the random collage models, called {\it dead leaves models} by the French school: see [L-M-H]. Here the $\phi_i$ are assumed to have bounded support, the terms have a random depth and, instead of being simply added, each term occludes anything behind it with respect to depth. This means $I(x,y)$ equals the one $\phi_i$ which is in front of all the others whose support contains $(x,y)$. This theory has major troubles with both infra-red and ultra-violet limits but it does provide the best approximation to date of the empirical statistics of images. It introduces explicitly the hidden variables describing the discrete objects in the image and allows one to model their preferred geometries.

Crafting models of this type is not simply mathematically
satisfying. It is central to the main application of computer
vision: object recognition. When an object of interest is obscured
in a cluttered badly lit scene, one needs a $p$-value for the
hypothesis test --- is this fragment of stuff part of the
sought-for object or an accidental conjunction of things occurring
in generic images? To get this $p$-value, one needs a null
hypothesis, a theory of generic images.

\subsection{Stochastic shapes via random diffeomorphisms and fluid flow}\vskip-5mm \hspace{5mm}

As we have seen in the last section, modeling images leads to
objects and these objects have shape --- so we need stochastic
models of shape, the ultimate non-linear sort of thing. Again it
is natural to consider this in the continuum limit and consider a
$k$-dimensional shape to be a subset of $\mathbb{R}^k$, e.g.\ a
connected open subset with nice boundary $\Gamma$. It is very
common in multiple images of objects like faces, animals, clothes,
organs in your body, to find not identical shapes but warped
versions. How is this to be modeled? One can follow the ideas of
the previous section and take a highly empirical approach,
gathering huge databases of faces or kidneys. This is probably the
road to the best pattern recognition in the long run. But another
principle that Grenander has always emphasized is to take
advantage of the group of symmetries of the situation --- in this
case, the group of all diffeomorphisms of $\mathbb{R}^k$. He and
Miller and collaborators (see [Gr-M]) were led to rediscover the
point of view of Arnold which we next describe.

Let ${\cal G}_n = \text{ group of diffeomorphisms on } \mathbb{R}^n$ and
${\cal SG}_n$ be the volume-preserving subgroup. We want to bypass issues of the exact degree of differentiability of these diffeomorphisms, but consider ${\cal G}_n$ and ${\cal SG}_n$ as infinite dimensional Riemannian manifolds. Let $\{\theta_t\}_{0\leq t\leq 1}$ be a path in ${\cal SG}_n$ and define its length by:
$$ \text{length of path } = \int \left( \sqrt{\int_{\mathbb{R}^n} \|  \frac{\partial \theta_t}{\partial t} (\theta^{-1}_t (x)) \|^2 d\vec{x}}\: \right) dt.$$
This length is nothing but the {\it right}-invariant Riemannian metric:
$$\text{dist}(\theta, (I+\epsilon \vec{v})\circ \theta)^2=\epsilon^2 \int \|\vec{v}\|^2 dx_1\cdot \cdot dx_n, \text{ where div}(\vec{v})\equiv 0.$$
Arnold's beautiful theorem is:

{\bf Theorem}  {\it Geodesics in ${\cal SG}_n$ are solutions of
Euler's equation:}
$$\frac{\partial v_t}{\partial t} +(v_t \cdot \nabla)v_t = \nabla p,\text{ some pressure } p.$$
This result suggests using geodesics on suitable infinite dimensional manifolds to model optimal warps between similar shapes in images and using diffusion on these manifolds to craft stochastic models. But we need to get rid of the volume-preserving restriction. The weak metric used by Arnold no longer works on the full ${\cal G}_n$ and in [C-R-M], Christensen et al introduced:
$$\|\vec{v}\|^2_L  = \int <L\vec{v} \cdot \vec{v}>dx_1 \cdots dx_n$$
where $v$ is any vector field and $L$ is a fixed positive self-adjoint differential operator e.g.\ $(I-\Delta)^m, m> n/2.$ Then a path $\{\theta_t\}$ in $G$ has both a velocity:
$$v_t=\frac{\partial \theta_t}{\partial t} (\theta^{-1}_t (x))$$
and a {\it momentum:} $u_t =Lv_t$ (so $v_t=K\ast u_t$, $K$ the Green's function of $L$). What is important here is that the momentum $u_t$ can be a generalized function, even when $v_t$ is smooth. The generalization of Arnold's theorem, first derived by Vishik, states that geodesics are:
$$\frac{\partial u_t}{\partial t} + (v_t \cdot \nabla) (u_t)+\text{div}(v_t) u_t =
-\sum_i(u_t)_i \vec{\nabla} ((v_t)_i).$$
This equation is a new kind of regularized compressible Euler equation, called by Marsden the template matching equation (TME). The left hand side is the derivative along the flow of the momentum, as a measure, and the right hand side is the force term.

A wonderful fact about this equation is that by making the momentum singular, we get very nice equations for geodesics on the ${\cal G}_n$-homogeneous spaces: \\
\hspace*{.3in}(a) ${\cal L}_n =$ set of all $N$-tuples of distinct points in $\mathbb{R}^n$ and\\
\hspace*{.3in}(b) ${\cal S}_n =$ set of all images of the unit ball under a diffeomorphism.\\
In the first case, we have ${\cal L}_n \cong {\cal G}_n/{\cal G}_{n,0}$
where ${\cal G}_{n,0}$ is the stabilizer of a specific set $\{P^{(0)}_1,\cdots, P^{(0)}_N\}$ of $N$ distinct points. To get geodesics on ${\cal L}_n$, we look for `particle solutions of the TME', i.e.\
$$\vec{u}_t=\sum^N_{i=1} \vec{u}_i(t) \delta_{P_i(t)}$$
where $\left\{P_1(t),\cdots, P_N(t)\right\}$ is a path in ${\cal L}_n$
The geodesics on ${\cal G}_n$, which are perpendicular to all cosets $\theta {\cal G}_{n,0}$, are then the geodesics on ${\cal L}_n$ for the quotient metric:
\begin{eqnarray*}
\text{dist}(\{P_i\}, \{P_i+\epsilon v_i\})^2 & =&
\epsilon^2 {\inf_{\genfrac{(}{)}{0pt}{}{v \text{ on }\mathbb{R}^n}{v(P_i)=v_i}} \int} <Lv,v>\\
& = &\epsilon^2 \displaystyle{\sum_{i,j}} G_{ij} (v_i\cdot v_j)
\end{eqnarray*}
where $G=K(\|P_i-P_j\|)^{-1}$. For these we get the interesting Hamiltonian ODE:
\begin{eqnarray*}
\frac{dP_i}{dt} &=& 2 \sum_j K(\|P_i-P_j\|)\vec{u}_j \\
\frac{du_i}{dt} &=& -\sum_j \nabla_{P_i} K(\|P_i-P_j\|)\cdot (\vec{u}_i \cdot \vec{u}_j)
\end{eqnarray*}
which makes points traveling in the same direction attract each other and points going in opposite directions repel each other. This space leads to a non-linear version of the theory of landmark points and shape statistics of Kendall [Sm] and has been developed by Younes [Yo].

A similar treatment can be made for the space of shapes ${\cal S}_n \cong {\cal G}_n/{\cal G}_{n,1}$, where ${\cal G}_{n,1}$ is the stabilizer of the unit sphere. Geodesics on ${\cal S}_n$ come from solutions of the TME for which $\vec{u}_t$ is supported on the boundary of the shape and perpendicular to it. Even though the first of these spaces ${\cal S}_2$ might seem to be quite a simple space, it seems to have a remarkable global geometry reflecting the many perceptual distinctions which we make when we recognize a similar shapes, e.g.\ a cell decomposition reflecting the different possible graphs which can occur as the `medial axis' of the shape. This is an area in which I anticipate interesting results. We can also use these Riemannian structures to define Brownian motion on ${\cal G}_n, {\cal S}_n$ and ${\cal L}_n$ (see [D-G-M], [Yi]). Putting a random stopping time on this walk, we get probability measures on these spaces. To make the ideas more concrete, in figure 6 we show a simulation of the random walk on ${\cal S}_2$.
\begin{figure}[htp]
\begin{center}
\epsfig{width=5in,file=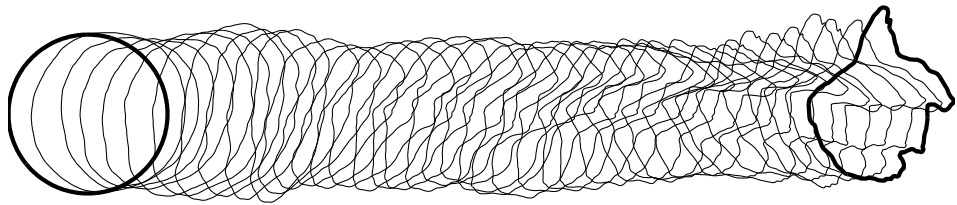}
%\center{\epsfig{width=5in,file=DiffEx3.eps}
\begin{minipage}[h]{10cm}
\caption{\footnotesize An example of a random walk in the space of 2D shapes ${\cal S}_2$. The initial point is
the circle on the left. A constant translation to the right has been added so the figures can be distinguished.
The operator $L$ defining the metric is $(I-\triangle)^2$}
\end{minipage}
\end{center}
\end{figure}

\vspace*{-1cm}

\section{Final thoughts}\vskip-5mm \hspace{5mm}

The patterns which occur in nature's sensory signals are complex but allow mathematical modeling. Their study has gone through several phases. At first, `off-the-shelf' classical models (e.g.\ linear Gaussian models) were adopted based only on intuition about the variability of the signals. Now, however, two things are happening: computers are large enough to allow massive data gathering to support fully non-parametric models. And the issues raised by these models are driving the study of new areas of mathematics and the development of new algorithms for working with these models. Applications like general purpose speech recognizers and computer driven vehicles are likely in the foreseeable future. Perhaps the ultimate dream is a fully unsupervised learning machine which is given only signals from the world and which finds their statistically significant patterns with no assistance: something like a baby in its first 6 months.

\label{lastpage}

\end{document}